\documentclass[reqno,12pt]{amsart}
\usepackage{amsmath,amsthm,amssymb,amsfonts,amscd}
\usepackage{xypic}
\theoremstyle{plain} 
	\newtheorem{thm}{Theorem}[section]
	\newtheorem*{thm*}{Theorem}
	
	\newtheorem{lem}[thm]{Lemma}
	\newtheorem{sublem}[thm]{Sub-Lemma}

	\newtheorem*{conj*}{Conjecture}
	
\theoremstyle{definition}

\theoremstyle{remark}
	\newtheorem{rem}[thm]{Remark}
	
	\newtheorem*{pf}{Proof}
\numberwithin{equation}{section}
\def\CC{{\mathbb C}}

\def\HH{{\mathbb H}}

\def\PP{{\mathbb P}}
\def\QQ{{\mathbb Q}}

\def\ZZ{{\mathbb Z}}

\def\I{{\mathcal I}}

\def\M{{\mathcal M}}

\def\X{{\mathcal X}}

\def\p{\partial }
\begin{document}
\title{Gromov--Witten invariants for mirror orbifolds 
of simple elliptic singularities}
\date{\today}
\author{Ikuo Satake}
\address{Faculty of Education, Kagawa University, 
1-1 Saiwai-cho Takamatsu Kagawa, 760-8522, Japan}
\email{satakeikuo@gmail.com}
\author{Atsushi Takahashi}
\address{Department of Mathematics, Graduate School of Science, Osaka University, 
Toyonaka Osaka, 560-0043, Japan}
\email{takahashi@math.sci.osaka-u.ac.jp}
\begin{abstract}
We consider a mirror symmetry of simple elliptic singularities. 
In particular, we construct isomorphisms of Frobenius manifolds among 
the one from the Gromov--Witten theory of a weighted projective line,
the one from the theory of primitive forms for a universal unfolding of a simple elliptic singularity 
and the one from the invariant theory for an elliptic Weyl group.
As a consequence, we give a geometric interpretation of the Fourier coefficients of 
an eta product considered by K. Saito.
\end{abstract}
\maketitle
\section*{Introduction}
Mirror symmetry can be understood as a duality between 
algebraic geometry and symplectic geometry.
It is an interesting problem to understand based on the philosophy of mirror symmetry 
some mysterious correspondences among
isolated singularities, root systems and discrete groups such as Schwartz's triangle groups.
Let $f(x,y,z)$ be a holomorphic function which has an isolated singularity 
only at the origin $0\in\CC^3$. 
A distinguished basis of vanishing cycles
in the Milnor fiber of $f$ can be categorified to an $A_\infty$-category ${\rm Fuk}^{\to}(f)$
called the directed Fukaya category whose derived category $D^b{\rm Fuk}^\to(f)$
is, as a triangulated category, an invariant of the holomorphic function $f$. 
If $f(x,y,z)$ is a weighted homogeneous polynomial 
then one can consider another interesting triangulated category, the category 
of a maximally-graded singularity $D^{L_f}_{Sg}(R_f)$:
\begin{equation}
D^{L_f}_{Sg}(R_f):=D^b({\rm gr}^{L_f}\text{-}R_f)/D^b({\rm proj}^{L_f}\text{-}R_f),
\end{equation}
where $R_f:=\CC[x,y,z]/(f)$ and $L_f$ is the maximal grading (see section one of \cite{et:1} for the definition) of $f$.
This category $D^{L_f}_{Sg}(R_f)$ is considered as an analogue 
of the bounded derived category of coherent sheaves on a smooth proper algebraic variety.
In this setting, homological mirror symmetry conjectures can be stated as follows:
\begin{conj*}[\cite{et:1}\cite{t:1}]
\begin{enumerate}
\item 
Let $f(x,y,z)$ be an invertible polynomial $($see section one of \cite{et:1} for the definition$)$.
There should exist a quiver with relations $(Q,I)$ and triangulated equivalences
\begin{equation}\label{hms:1}
D^{L_f}_{Sg}(R_f)\simeq D^b({\rm mod}\text{-}\CC Q/I)\simeq D^b{\rm Fuk}^\to(f^t),
\end{equation}
where $f^t$ denotes the Berglund--H\"{u}bsch transpose of $f$.
\item There should exist triangulated equivalences
\begin{equation}\label{hms:2}
D^b{\rm coh}(\PP^1_{a_1,a_2,a_3})\simeq D^b({\rm mod}\text{-}\CC Q_{a_1,a_2,a_3}/I')\simeq
D^b{\rm Fuk}^\to (T_{a_1,a_2,a_3}),
\end{equation}
where $\PP^1_{a_1,a_2,a_3}$ is the orbifold $\PP^1$ with $3$ isotropic points of orders 
$a_1,a_2,a_3$, $Q_{a_1,a_2,a_3}$ is a quiver given by the following graph 
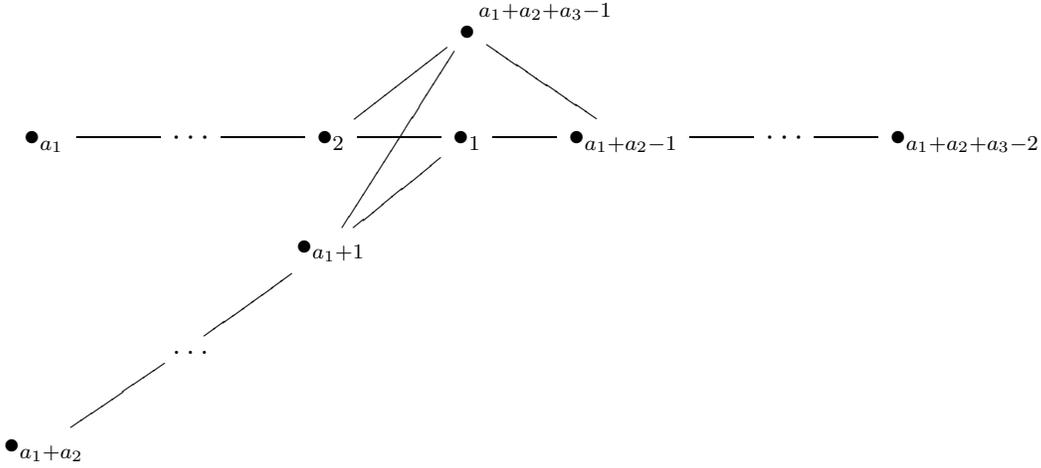
\[
\xymatrix{ 
 & & & {\bullet}  \ar@{-}[dr]  \ar@{-}[ldd] \ar@{}^{a_1+a_2+a_3-1}[r] 
 & & &  \\
 {\bullet}_{a_1} \ar@{-}[r]  & {\cdots} \ar@{-}[r]  & {\bullet}_2 \ar@{-}[r] \ar@{-}[ur] & {\bullet}_1 \ar@{-}[dl] \ar@{-}[r] \ar@{-}[r] & 
{\bullet}_{a_1+a_2-1} \ar@{-}[r]  & {\cdots} \ar@{-}[r]  &{\bullet}_{a_1+a_2+a_3-2}   \\
& &  {\bullet}_{a_1+1} \ar@{-}[dl] & & & &  \\
 & {\dots} \ar@{-}[dl] & & & & & \\
{\bullet}_{a_1+a_2}& & & & & &
  }
\]
with the orientation from vertices with smaller indices to those with larger indices 
and $I'$ is the ideal generated by two generic paths from the $1$-st vertex to 
the $a_1+a_2+a_3-1$-th vertex,
and 
$T_{a_1,a_2,a_3}:=x_1^{a_1}+x_2^{a_2}+x_3^{a_3}-cx_1x_2x_3$, $c\in\CC^*$.
\end{enumerate}
\end{conj*}
It is natural to expect the following from (ii) of the above homological mirror symmetry conjectures 
since their ``complexified K\"ahler moduli spaces" should be isomorphic and 
there should exist Frobenius structures (K.~Saito's flat structures) on them:
\begin{conj*}
There should exist isomorphisms of Frobenius manifolds $($see for example \cite{d:1,st:1} for the definition$)$ among 
\begin{enumerate}
\item $M_{\PP^1_{a_1,a_2,a_3}}$, 
the one constructed from the Gromov--Witten theory of $\PP^1_{a_1,a_2,a_3}$,
\item $M_{(Q_{a_1,a_2,a_3}, I')}$, 
the one constructed from the invariant theory of the reflection group 
associated to the quiver with relations $(Q_{a_1,a_2,a_3}, I')$,
\item $M_{T_{a_1,a_2,a_3},\infty}$,
the one constructed from the universal unfolding of 
$T_{a_1,a_2,a_3}$ by the choice of primitive form ``at $c= \infty$''.
\end{enumerate}
\end{conj*}
\begin{rem}
It is also a part of conjecture that there exist Frobenius manifolds $M_{(Q_{a_1,a_2,a_3}, I')}$ for 
$1/a_1+1/a_2+1/a_3<1$.
\end{rem}
Rossi shows in \cite{r:1} that Conjecture holds under the condition $1/a_1+1/a_2+1/a_3>1$.
The next case to consider is when the triple $(a_1,a_2,a_3)$ satisfies 
the condition $1/a_1+1/a_2+1/a_3=1$, in other words, the case when the polynomial 
$f$ defines a simple elliptic singularity (see \cite{et:1} for this relation between 
$(a_1,a_2,a_3)$ and $f$). 
In particular, in this paper we shall give a proof of 
the above Conjecture for $(a_1,a_2,a_3)=(3,3,3)$ 
with the explicit presentation of the potential which 
gives us interesting quasi-modular forms 
based on the uniqueness of 
the solution of the WDVV equation.
The following is our main result in this paper:
\begin{thm*}
We have isomorphisms of Frobenius manifolds
\[
M_{\PP^1_{3,3,3}}\simeq M_{E_6^{(1,1)}}\simeq M_{T_{3,3,3},\infty},
\]
where $M_{E_6^{(1,1)}}$ denotes the Frobenius manifold 
constructed from the invariant theory of the elliptic Weyl group of type $E_6^{(1,1)}$. 
Moreover, the genus zero Gromov--Witten potential $F_0^{\PP^1_{3,3,3}}$ and 
the genus one Gromov--Witten potential $F_1^{\PP^1_{3,3,3}}$, which is also 
considered as the $G$-function $($see \cite{dz:2} for the definition$)$ on $M_{E_6^{(1,1)}}$ and as the one on $M_{T_{3,3,3},\infty}$, 
are expressed by quasi-modular forms. 
\qed
\end{thm*}
An important consequence of this theorem is that 
we can give a geometric interpretation of the Fourier coefficients of 
an eta product considered by K. Saito \cite{sa:eta}:
\begin{thm*}
Denote by $\eta(\tau)$ the Dedekind's eta function 
\[
\eta(\tau):=e^{\frac{2\pi\sqrt{-1}\tau}{24}}\prod_{n\ge 1}
\left(1-e^{2\pi\sqrt{-1}n\tau}\right),\quad \tau\in \HH:=\{\tau\in\CC~|~{\rm Im}\tau >0\}.
\]
The eta product $\eta(3\tau)^3/\eta(\tau)$ is a generating function of Gromov--Witten 
invariants of $\PP^1_{3,3,3}$.
More precisely, the Fourier coefficient $c_k$ defined by 
\begin{equation}
\frac{\eta(3\tau)^3}{\eta(\tau)}=e^{\frac{2\pi\sqrt{-1}\tau}{3}}\sum_{k\ge 0}c_k 
e^{2\pi\sqrt{-1}k\tau}
\end{equation}
is the Gromov--Witten invariant 
\[
\int_{[\overline{\M}_{0,0,3k[\PP^1_{3,3,3}]}(\PP^1_{3,3,3})]^{vir}}
ev_1^*\gamma_1\wedge ev_2^*\gamma_2\wedge
ev_3^*\gamma_3,
\]
where $\gamma_i$ is an element of $H^{2/3}_{orb}(\PP^1_{3,3,3},\QQ)$
corresponding to the $i$-th isotropic point on $\PP^1_{3,3,3}$.
\qed
\end{thm*}
We can also apply the same method to prove the Conjecture for the two other cases when $(a_1,a_2,a_3)=(2,4,4), (2,3,6)$. 
However, we omit them here since 
the number of monomials in those potentials are large 
(more than 50 for $(2,4,4)$ and more than 200 for $(2,3,6)$) 
we can not give the explicit presentation of the potential in this paper 
and we could understand not all but a few of interesting quasi-modular forms 
appearing in those potentials. 
We can also consider a similar problem for which we do not have a 
hypersurface singularity:
\begin{thm*}
We have an isomorphism of Frobenius manifolds
\[
M_{\PP^1_{2,2,2,2}}\simeq M_{D_4^{(1,1)}},
\]
where $\PP^1_{2,2,2,2}$ denotes an orbifold $\PP^1$ with four isotropic points 
of orders $2$ and $M_{D_4^{(1,1)}}$ denotes the Frobenius manifold 
constructed from the invariant theory of the elliptic Weyl group of type $D_4^{(1,1)}$. 
Moreover, the genus zero Gromov--Witten potential $F_0^{\PP^1_{2,2,2,2}}$ and 
the genus one Gromov--Witten potential $F_1^{\PP^1_{2,2,2,2}}$, which is also 
considered as the $G$-function on $M_{D_4^{(1,1)}}$,
are expressed by quasi-modular forms. 
\qed
\end{thm*}
Note that in order to obtain the mirror isomorphism 
we have to develop the theory of primitive forms for 
a pair consisting in a singularity and its symmetry group.
Once we have such a theory, we may apply it for the pair $(T_{2,4,4},\ZZ/2\ZZ)$, for example. 
If the triple $(a_1,a_2,a_3)$ satisfies the condition $1/a_1+1/a_2+1/a_3=1$, 
then we have the triangulated equivalence 
$D^{L_f}_{Sg}(R_f)\simeq D^b{\rm coh}(\PP^1_{a_1,a_2,a_3})$
of Buchweitz--Orlov type (see \cite{u:1}).
Also note that a mathematical formulation of 
the topological A-model for Landau--Ginzburg orbifold theory 
is considered in \cite{fjr:1}, which is called Fan--Jarvis--Ruan--Witten (FJRW) theory. 
Therefore, it is also natural to consider the following:
\begin{conj*}
Let $T_{a_1,a_2,a_3}$ be a polynomial which defines a simple elliptic singularity 
$\widetilde{E}_6$, $\widetilde{E}_7$ or $\widetilde{E}_8$. 
There should exist an isomorphism of Frobenius manifolds between
\begin{enumerate}
\item $M_{(T_{a_1,a_2,a_3},\ZZ/d\ZZ), FJRW}$, 
the one constructed from the FJRW theory for the pair $(T_{a_1,a_2,a_3},\ZZ/d\ZZ)$ 
where $d=6,7,8$ for $\widetilde{E}_6,\widetilde{E}_7,\widetilde{E}_8$ respectively,
\item $M_{T_{a_1,a_2,a_3},0}$,
the one constructed from the universal unfolding of 
$T_{a_1,a_2,a_3}$ by the choice of primitive form ``at $c=0$''.
\end{enumerate}
\end{conj*}
The authors are notified that 
Krawitz--Shen \cite{ks:1} gives a proof of this Conjecture based on the calculations of $M_{T_{a_1,a_2,a_3},0}$
by Noumi--Yamada \cite{ny:1} and 
Milanov--Ruan \cite{mr:1} prove a generalization of this, namely, the one for all genus potentials 
and their quasi-modularity.
\bigskip
\noindent
{\bf Acknowledgement}\\
\indent
The second named author is supported 
by JSPS KAKENHI Grant Number 20360043, 24684005. 
We thank the anonymous referee for carefully reading our paper.
\section{Gromov--Witten theory for orbifolds}
Gromov--Witten theory is generalized for orbifolds (smooth proper Deligne--Mumford stacks). 
It is first studied by Chen--Ruan \cite{cr:1} in symplectic geometry and 
later by Abramovich-Graber--Vistoli \cite{agv:1} in algebraic geometry.
In order to generalize Gromov--Witten theory for manifolds to the one for orbifolds,
one also needs to count the number of ``stable maps from orbifold curves".
For this purpose, in \cite{cr:1} the notion of orbifold stable maps is introduced 
and in \cite{agv:1} the notion of twisted stable maps is introduced.
These two constructions are quite different, however, as the usual Gromov--Witten 
theory for manifolds,  
they are expected to give the same Gromov--Witten invariants since they have common philosophy.
In this paper, we will introduce  Gromov--Witten invariants following \cite{agv:1} for simplicity. 
Let $\X$ be an orbifold $($or a smooth proper Deligne--Mumford stack over $\CC)$.
Then, for $g\in\ZZ_{\ge 0}$, $n\in\ZZ_{\ge 0}$ 
and $\beta\in H_2(\X,\ZZ)$, 
the moduli space (stack) $\overline{\M}_{g,n,\beta}(\X)$ 
of orbifold (twisted) stable maps of genus $g$ 
with $n$-marked points of degree $\beta$ is defined.
There exists a virtual fundamental class $[\overline{\M}_{g,n,\beta}(\X)]^{vir}$ 
and Gromov--Witten invariants of genus $g$ 
with $n$-marked points of degree $\beta$ are defined as usual except for 
that we have to use the orbifold cohomology group $H^*_{orb}(\X,\QQ)$:
\[
\left<\gamma_1,\dots, \gamma_n\right>_{g,n,\beta}^\X:=
\int_{[\overline{\M}_{g,n,\beta}(\X)]^{vir}}ev_1^*\gamma_1\wedge \dots 
ev_n^*\gamma_n,\quad \gamma_1,\dots,\gamma_n\in H^*_{orb}(\X,\QQ),
\]
where $ev^*_i:H^*_{orb}(\X,\QQ)\longrightarrow H^*(\overline{\M}_{g,n,\beta}(\X),\QQ)$
denotes the induced homomorphism by the evaluation map.
We also consider the generating function 
\[
F_g^\X:=\sum_{n,\beta}\frac{1}{n!}\left<{\bf t},\dots, {\bf t}\right>_{g,n,\beta}^\X q^\beta,
\quad {\bf t}=\sum_{i}t_i\gamma_i
\]
and call it the genus $g$ potential where $\{\gamma_i\}$ denotes a $\QQ$-basis of 
$H^*_{orb}(\X,\QQ)$.
The main result in \cite{agv:1} and \cite{cr:1} tell us that 
we can treat the Gromov--Witten theory defined for orbifolds 
as if $\X$ were usual manifold. 
In particular, we have the point axiom,
the divisor axiom for a class in $H^2(\X,\QQ)$
and the associativity of the quantum product, namely, the WDVV equation 
(see, for example, \cite{cr:1} for details of these axioms.), 
which gives a (formal) Frobenius manifold.
These axioms enable us to calculate genus zero Gromov--Witten potential $F_0^\X$ easily.
In this paper, we shall only consider the case when $\X$ is $\PP^1_{2,2,2,2}$
or $\PP^1_{3,3,3}$, the orbifold $\PP^1$ with $4$ isotropic points of order $2$
or the orbifold $\PP^1$ with $3$ isotropic points of order $3$.
Note that both are given by the global quotient of an elliptic curve ${\bf E}$,
more precisely, we have $\PP^1_{2,2,2,2}=[{\bf E}/(\ZZ/2\ZZ)]$ and 
$\PP^1_{3,3,3}=[{\bf E}/(\ZZ/3\ZZ)]$.
For these examples, 
by the uniqueness result on genus zero and one potentials,
we shall see that the two definitions of Gromov--Witten invariants 
by \cite{agv:1} and \cite{cr:1} coincides. 
\section{Explicit calculations for $\PP^1_{2,2,2,2}$}
The orbifold cohomology group of $\PP^1_{2,2,2,2}$ is, as a vector space, just the
singular cohomology group of the inertia orbifold
\[
\I\PP^1_{2,2,2,2}=\PP^1_{2,2,2,2} \bigsqcup B(\ZZ/2 \ZZ)\bigsqcup B(\ZZ/2 \ZZ)\bigsqcup B(\ZZ/2 \ZZ)\bigsqcup B(\ZZ/2 \ZZ),
\]
and the orbifold Poincar\'{e} pairing is given by twisting the usual Poincar\'{e} pairing:
\[
\displaystyle
\int_{\PP^1_{2,2,2,2}} \alpha \cup_{orb} \beta := \int_{\I\PP^1_{2,2,2,2}} \alpha \cup I \beta,
\]
where $I$ is the involution defined in \cite{agv:1, cr:1}.
Therefore, we can choose a basis $\gamma_0,\dots, \gamma_5$ of the orbifold cohomology group 
$H^*_{orb}(\PP^1_{2,2,2,2},\QQ)$ such that 
\[
H^{0}_{orb}(\PP^1_{2,2,2,2},\QQ) \simeq \QQ\gamma_0,\quad 
H^{1}_{orb}(\PP^1_{2,2,2,2},\QQ) \simeq \bigoplus_{i=1}^4\QQ\gamma_i,\quad 
H^{2}_{orb}(\PP^1_{2,2,2,2},\QQ) \simeq \QQ\gamma_5,
\]
and
\[
\int_{\PP^1_{2,2,2,2}}\gamma_0\cup\gamma_5=1,\quad 
\int_{\PP^1_{2,2,2,2}}\gamma_i\cup\gamma_j=\frac{1}{2}\delta_{i,j},\ i,j=1,\dots 4.
\]
Denote by $t_0,\dots, t_5$ the dual coordinates of the $\QQ$-basis $\gamma_0,\dots, \gamma_5$.
In discussion below, by applying the divisor axiom, 
we consider $\log q$ as a flat coordinate instead of $t_5$. 
\subsection{Genus zero potential}
\begin{thm}\label{S1130-1}
The genus zero Gromov--Witten potential $F_0^{\PP^1_{2,2,2,2}}$ of $\PP^1_{2,2,2,2}$ 
is given as follows$:$
\begin{align*}
F_0^{\PP^1_{2,2,2,2}}=&\frac{1}{2}t_0^2 \log q+\frac{1}{4}t_0 (t_1^2+t_2^2+t_3^2+t_4^2)\\
&+(t_1 t_2 t_3 t_4)\cdot f_0(q)+\frac{1}{4}(t_1^4+t_2^4+t_3^4+t_4^4)\cdot f_1(q)\\
&+\frac{1}{6}(t_1^2 t_2^2+t_1^2 t_3^2+t_1^2 t_4^2+t_2^2 t_3^2
+t_2^2 t_4^2+t_3^2 t_4^2)\cdot f_2(q),
\end{align*}
where 
\begin{align}
&f_0(q):= \frac{1}{2}\left(f(q)-f(-q)\right),\\
&f_1(q):= f(q^4),\\
&f_2(q):= f(q)-f_0(q)-f_1(q),\\
&f(q):=-\frac{1}{24}+\sum_{n=1}^\infty n\frac{q^n}{1-q^n}
=-q\frac{d}{dq}\log(\eta(q)),\\
&\eta(q):=q^{\frac{1}{24}}\prod_{n=1}^{\infty}(1-q^n).
\end{align}
\end{thm}
\begin{pf}
We can deduce Theorem~\ref{S1130-1}
from the following uniqueness property of the potential:
\begin{lem}\label{S1130-2}
There exists a unique $6$-dimensional formal Frobenius structure with 
flat coordinates $t_0,t_1,t_2,t_3,t_4,t$ satisfying the following conditions$:$
\begin{enumerate}
\item The Euler vector field $E$ is given by  
$E=t_0\frac{\p}{\p t_0}+\sum_{k=1}^4\frac{1}{2}t_k\frac{\p}{\p t_k}$.
\item The Frobenius potential $F_0$ is given by
\begin{align*}
F_0=&\frac{1}{2}t_0^2 t+\frac{1}{4}t_0 (t_1^2+t_2^2+t_3^2+t_4^2)\\
&+(t_1 t_2 t_3 t_4)\cdot f_0(e^t)+\frac{1}{4}(t_1^4+t_2^4+t_3^4+t_4^4)\cdot f_1(e^t)\\
&+\frac{1}{6}(t_1^2 t_2^2+t_1^2 t_3^2+t_1^2 t_4^2+t_2^2 t_3^2
+t_2^2 t_4^2+t_3^2 t_4^2)\cdot f_2(e^t),
\end{align*}
where $f_0(q),f_1(q),f_2(q)$ have the following formal power series expansions$:$
\[
f_0(q)=\sum_{n=1}^\infty a_nq^{n}\text{ with }a_1=1,\quad 
f_1(q)=\sum_{n=0}^\infty b_nq^{n},\quad 
f_2(q)=\sum_{n=0}^\infty c_nq^{n}.
\]
\end{enumerate}
\end{lem}
\begin{pf}
We can show that the WDVV equation is equivalent to the following differential equations$:$ 
\begin{align}
&q\frac{d}{dq}f_0(q)=\frac{8}{3}f_0(q)f_2(q)-24f_0(q)f_1(q),\\
&q\frac{d}{dq}f_1(q)=-\frac{2}{3}f_0(q)^2-\frac{16}{3} f_1(q) f_2(q)+\frac{8}{9}f_2(q)^2,\\
&q\frac{d}{dq}f_2(q)=6f_0(q)^2-\frac{8}{3}f_2(q)^2.
\end{align}
Hence, we have the following recursion relations for $a_n,b_n,c_n$:
\begin{align}
&na_n=\frac{8}{3}\sum_{k=1}^na_kc_{n-k}-24\sum_{k=1}^na_kb_{n-k},\\
&nb_n=-\frac{2}{3}\sum_{k=1}^{n-1}a_ka_{n-k}-\frac{16}{3}\sum_{k=0}^nb_kc_{n-k}
+\frac{8}{9}\sum_{k=0}^nc_kc_{n-k},\\
&nc_n=6\sum_{k=1}^{n-1}a_k a_{n-k}-\frac{8}{3}\sum_{k=0}^n c_k c_{n-k}.
\end{align}
In particular, by setting $n=0,1$, we get $c_0=0$ and $b_0=-1/24$. 
Therefore, the above recursion relations have the unique solution.
\qed
\end{pf}
Next, we construct the analytic solution to the WDVV equation as follows.
\begin{lem}\label{22}
Put 
\begin{align}
&f_0(q):= \frac{1}{2}\left(f(q)-f(-q)\right),\\
&f_1(q):= f(q^4),\\
&f_2(q):= f(q)-f_0(q)-f_1(q),\\
&f(q):=-\frac{1}{24}+\sum_{n=1}^\infty n\frac{q^n}{1-q^n}
=-q\frac{d}{dq}\log(\eta(q)).
\end{align}
Then the functions $f_0(q),\,f_1(q),\,f_2(q)$ 
satisfies the following differential equations$:$
\begin{align}
&q\frac{d}{dq}f_0(q)=\frac{8}{3}f_0(q)f_2(q)-24f_0(q)f_1(q),\\
&q\frac{d}{dq}f_1(q)=-\frac{2}{3}f_0(q)^2-\frac{16}{3} f_1(q) f_2(q)+\frac{8}{9}f_2(q)^2,\\
&q\frac{d}{dq}f_2(q)=6f_0(q)^2-\frac{8}{3}f_2(q)^2.
\end{align}
\end{lem}
\begin{pf}
Put 
\begin{align*}
&\vartheta_2(q):=\sum_{m \in \ZZ}q^{(m+\frac{1}{2})^2},\\
&\vartheta_3(q):=\sum_{m \in \ZZ}q^{m^2},\\
&\vartheta_4(q):=\sum_{m \in \ZZ}(-1)^m q^{m^2},\\
&X_i(q):=q\frac{d}{dq}\log \vartheta_i\ (i=2,3,4).
\end{align*}
Then the following differential relations
\begin{align}
&\frac{1}{2}q\frac{d}{dq}(X_2(q)+X_3(q))=2X_2(q) X_3(q),\\
&\frac{1}{2}q\frac{d}{dq}(X_3(q)+X_4(q))=2X_3(q) X_4(q),\\
&\frac{1}{2}q\frac{d}{dq}(X_4(q)+X_2(q))=2X_4(q) X_2(q)
\end{align}
are classically known as Halphen's equations (see \cite{S1130-Ohyama1}). 
For the proof of Lemma~\ref{22}, we should only prove that 
\begin{align}
&X_2(q)=-6 f_1(q)+\frac{2}{3}f_2(q),
\label{S1129-8}\\
&X_3(q)=2 f_0(q)-\frac{4}{3}f_2(q),
\label{S1129-9}\\
&X_4(q)=-2 f_0(q)-\frac{4}{3}f_2(q).
\label{S1129-10}
\end{align}
We have
\begin{align}
&X_2(q)=q\frac{d}{dq}
\log[2\eta(q^2)^{-1}\eta(q^4)^{2}],
\label{S1129-2}\\
&X_3(q)=q\frac{d}{dq}
\log[\eta(q)^{-2}\eta(q^2)^{5}\eta(q^4)^{-2}],
\label{S1129-3}\\
&X_4(q)=q\frac{d}{dq}
\log[\eta(q)^2\eta(q^2)^{-1}]
\label{S1129-4}
\end{align}
by Jacobi's triple product formula (see \cite{Mumford}). 

For $f_0(q),\,f_1(q),\,f_2(q)$, we prepare the following 
Sub-Lemma.
\begin{sublem}\label{sublem:E2 identity}
For $f(q)$, 
we have 
\begin{equation}
\frac{1}{2}(f(q)+f(-q))=3 f(q^2)-2 f(q^4).
\label{S1129-1}
\end{equation}
\end{sublem}
\begin{pf}
We define $\sigma(n)\,(n\geq 1)$ by 
\[
f(q)=-\frac{1}{24}+\sum_{n=1}^{\infty}\sigma(n)q^n.
\]
If $n=2^k m$ with $m$ odd, then 
$\sigma(n)=(1+2+\cdots+2^k)\sigma(m)
=(2^{k+1}-1)\sigma(m)$. 
Thus we have $\sigma(2n)=3\sigma(n)$ if $n$ is odd. 
Also if $n$ is general, we have 
$\sigma(4n)=3\sigma(2n)-2\sigma(n)$.
Then we have (\ref{S1129-1}).
\qed
\end{pf}
By (\ref{S1129-1}), we have
\begin{align}
&f_0(q)=-q\frac{d}{dq}
\log[\eta(q)\eta(q^2)^{-\frac{3}{2}}\eta(q^4)^{\frac{1}{2}}],
\label{S1129-5}\\
&f_1(q)=-q\frac{d}{dq}
\log[\eta(q^4)^{\frac{1}{4}}],
\label{S1129-6}\\
&f_2(q)=-q\frac{d}{dq}
\log[\eta(q^2)^{\frac{3}{2}}\eta(q^4)^{-\frac{3}{4}}].
\label{S1129-7}
\end{align}
From (\ref{S1129-2})--(\ref{S1129-4}) and (\ref{S1129-5})--(\ref{S1129-7}), 
we have (\ref{S1129-8})--(\ref{S1129-10}).
\qed
\end{pf}

It is easy to show that, by our choice of basis $\gamma_0,\dots, \gamma_5$ of $H^*_{orb}(\PP^1_{2,2,2,2},\CC)$ 
and their dual coordinates $t_0,\dots, t_5=\log q$ in the beginning of this section,  
the Gromov--Witten potential is of the form in Lemma~\ref{S1130-2} except for the condition $a_1=1$.
The condition $a_1=1$ follows from the fact that 
the Gromov--Witten invariant $a_1$ counts the number of  
morphisms from $\PP^1_{2,2,2,2}$ to $\PP^1_{2,2,2,2}$ of degree one,
which is exactly the identity map. Hence, we have $a_1=1$.
Now, the statement in Theorem \ref{S1130-1} follows from the uniquness of the potential.
\qed
\end{pf}

By Theorem \ref{S1130-1}, the 
Gromov--Witten potential $F_0^{\PP^1_{2,2,2,2}}$ converges 
on the domain $|q|<1$. 
Thus it gives a Frobenius manifold 
$M_{\PP^1_{2,2,2,2}}\simeq 
\{z \in \CC\,|\,\text{Re}z<0\,\}\times \CC^5$
with flat coordinates $(\log q, t_0,t_1,t_2,t_3,t_4)$. 

For the elliptic root system of type $D_4^{(1,1)}$
(\cite{extendedII}), 
the domain ${\mathbb E}_{D_4^{(1,1)}}$ and 
the elliptic Weyl group $W_{D_4^{(1,1)}}$ are defined and 
the quotient space $M_{D_4^{(1,1)}}:={\mathbb E}_{D_4^{(1,1)}}//W_{D_4^{(1,1)}}
\simeq 
\{z \in \CC\,|\,\text{Re}z<0\,\}\times \CC^5$
has a structure of the Frobenius manifold 
(\cite{extendedII}, \cite{SatakeFrobenius}). 
Its potential is explicitly calculated in \cite{D4} as follows$:$
\begin{lem}{\rm{(\cite{D4})}}
By choosing the flat coordinates 
$t,e_0,e_1,e_3,e_4,e_2$ of 
$M_{D_4^{(1,1)}}$, the potential $F_0^{D_4^{(1,1)}}$ is expressed as 
\begin{align*}
F_0^{D_4^{(1,1)}}=&{1 \over 2}t (e_2)^2 \cr
&+{1 \over 4}e_2 [e_0^2+e_1^2+e_3^2+e_4^2]\cr
&+(e_0e_1e_3e_4)\cdot h_0(t)\\
&+{1 \over 4}(e_0^4+e_1^4+e_3^4+e_4^4)\cdot h_1(t)\cr
&+{1 \over 6}
(e_0^2e_1^2+e_0^2e_3^2+e_0^2e_4^2+e_1^2e_3^2 +e_1^2e_4^2+e_3^2e_4^2)
\cdot h_2(t),
\end{align*}
where 
\begin{align*}
&h_0(t)=\frac{1}{8}\Theta_{\omega_1,1}(e^t),\\
&h_1(t)=-\frac{1}{2}
\left[
\frac{1}{2}\frac{\frac{d}{dt}[\eta(e^{2t})]}{\eta(e^{2t})}
+{1 \over 24}\Theta_{0,1}(e^{t})
\right],\\
&h_2(t)=-\frac{3}{2}
\left[
\frac{1}{2}\frac{\frac{d}{dt}[\eta(e^{2t})]}{\eta(e^{2t})}
-{1 \over 24}\Theta_{0,1}(e^{t})
\right],\\
&\Theta_{0,1}(q)=\sum_{\gamma \in M}q^{(\gamma,\gamma)}=1+\cdots,\\
&\Theta_{\omega_1,1}(q)=
\sum_{\gamma \in M+\omega_1}q^{(\gamma,\gamma)}=8q+\cdots,
\end{align*}
where $M$ is the coroot lattice of $D_4$ and 
$\omega_1$ is the first fundamental weights in the notation of Bourbaki. 
\qed
\end{lem}

\begin{rem}
We remark that 
the correspondence of the above coordinates 
with the ones in \cite{D4} is 
\[
t=\pi\sqrt{-1}\tau,
e_0=c_0,e_1=c_1,e_3=c_3,e_4=c_4,e_2=\frac{-1}{2(2\pi\sqrt{-1})^2}c_2
\]
and we take the intersection form of the Frobenius manifold 
as $\frac{-1}{(2\pi\sqrt{-1})^2}I^*$ instead of $I^*$. 
\end{rem}

Since the potential $F_0^{D_4^{(1,1)}}$ satisfies the assumptions 
of the Lemma \ref{S1130-2}, we have 
\begin{thm}
The Frobenius manifold $M_{\PP^1_{2,2,2,2}}$ 
and the Frobenius manifold $M_{D_4^{(1,1)}}$ are 
isomorphic as Frobenius manifolds. 
\qed
\end{thm}
\subsection{Genus one potential}
We shall also give the genus one Gromov--Witten potential.
\begin{thm}\label{thm:D4-g1}
The genus one Gromov--Witten potential $F_1^{\PP^1_{2,2,2,2}}$ of $\PP^1_{2,2,2,2}$ 
is given as
\begin{equation}
F_1^{\PP^1_{2,2,2,2}}=-\frac{1}{2}\log(\eta(q^2)).
\end{equation}
\end{thm}
\begin{pf}
The first derivative of the genus one Gromov--Witten potential 
$q\frac{d}{dq}F_1^{\PP^1_{2,2,2,2}}$ is an element of $\QQ[[q]]$
since the Euler vector field is given by 
$E=t_0\frac{\p}{\p t_0}+\sum_{k=1}^4\frac{1}{2}t_k\frac{\p}{\p t_k}$, 
 $E\left(q\frac{d}{dq}F_1^{\PP^1_{2,2,2,2}}\right)=0$ 
and we have the divisor axiom.
Therefore, we only have to consider the (orbifold) stable maps 
with one marked point from 
smooth elliptic curves to $\PP^1_{2,2,2,2}=[{\bf E}/(\ZZ/2\ZZ)]$, which 
factor through the elliptic curve ${\bf E}$ by definition.
In particular, the number of coverings of degree $n$ from an elliptic curve to ${\bf E}$ 
is given by $\sigma(n):=\sum_{k|n}k$.
Hence, we have
\[
q\frac{d}{d q}F_1^{\PP^1_{2,2,2,2}}=f(q^2)=
-\frac{1}{24}+\sum_{n=1}^{\infty}\sigma(n)q^{2n}.
\]
One may also obtain the statement by Dubrovin--Zhang's Virasoro constraint \cite{dz:1}.
Indeed, Proposition 4 in \cite{dz:1} gives us the equation
\[
q\frac{d}{d q}F_1^{\PP^1_{2,2,2,2}}=f_1(q)+\frac{1}{3}f_2(q).
\]
By Sub-Lemma \ref{sublem:E2 identity}, we have $f_1(q)+\frac{1}{3}f_2(q)=f(q^2)$.
\qed
\end{pf}
The proof of Theorem \ref{thm:D4-g1} also shows that the genus one potential 
is uniquely reconstructed from the genus zero potential. 
In particular, this implies the $G$-function of $M_{D_4^{(1,1)}}$ coincides with 
$F_1^{\PP^1_{2,2,2,2}}$.
\section{Explicit calculations for $\PP^1_{3,3,3}$}
The orbifold cohomology group of $\PP^1_{3,3,3}$ is, as a vector space, just the
singular cohomology group of the inertia orbifold
\[
\I\PP^1_{3,3,3}=\PP^1_{3,3,3} \bigsqcup B(\ZZ/3 \ZZ)\bigsqcup B(\ZZ/3 \ZZ)\bigsqcup B(\ZZ/3 \ZZ),
\]
and the orbifold Poincar\'{e} pairing is given by twisting the usual Poincar\'{e} pairing:
\[
\displaystyle
\int_{\PP^1_{3,3,3}} \alpha \cup_{orb} \beta := \int_{\I\PP^1_{3,3,3}} \alpha \cup I \beta,
\]
where $I$ is the involution defined in \cite{agv:1, cr:1}.
Therefore, we can choose a $\QQ$-basis $\gamma_0,\dots, \gamma_7$ of the orbifold cohomology group 
$H^*_{orb}(\PP^1_{3,3,3},\QQ)$ such that 
\[
H^{0}_{orb}(\PP^1_{3,3,3},\QQ) \simeq \QQ\gamma_0,\quad 
H^{\frac{2}{3}}_{orb}(\PP^1_{3,3,3},\QQ) \simeq \bigoplus_{i=1}^3\QQ\gamma_i,
\]
\[ 
H^{\frac{4}{3}}_{orb}(\PP^1_{3,3,3},\QQ) \simeq \bigoplus_{i=4}^6\QQ\gamma_i,\quad 
H^{2}_{orb}(\PP^1_{3,3,3},\QQ) \simeq \QQ\gamma_7,
\]
and
\[
\int_{\PP^1_{3,3,3}}\gamma_0\cup\gamma_7=1,\quad 
\int_{\PP^1_{3,3,3}}\gamma_i\cup\gamma_j=\frac{1}{3}\delta_{i+j-7,0},\ i,j=1,\dots, 6.
\]
Denote by $t_0,\dots, t_7$ the dual coordinates of the $\QQ$-basis $\gamma_0,\dots, \gamma_7$.
In the discussion below, by applying the divisor axiom, 
we consider $\log q$ as a flat coordinate instead of $t_7$. 
\subsection{Genus zero potential}
\begin{thm}\label{thm:genus zero E6}
The genus zero Gromov--Witten potential $F_0^{\PP^1_{3,3,3}}$ of $\PP^1_{3,3,3}$ 
is given as follows$:$
\begin{align*}
F_0^{\PP^1_{3,3,3}}
=&\frac{1}{2}t_0^2 \log q
+\frac{1}{3}t_0 (t_1 t_6+t_2 t_5+t_3 t_4)+(t_1 t_2 t_3)\cdot f_0(q)\\
&+\frac{1}{6}(t_1^3+t_2^3+t_3^3)\cdot f_1(q)
+(t_1t_2t_5t_6+t_1t_3t_4t_6+t_2 t_3t_4t_5)\cdot f_2(q)\\
&+\frac{1}{2}(t_1^2t_4t_5+t_2^2t_4t_6+t_3^2t_5t_6)\cdot f_3(q)\\
&+\frac{1}{2}(t_1t_2t_4^2+t_1t_3t_5^2+t_2t_3t_6^2)\cdot f_4(q)
+\frac{1}{4}(t_1^2t_6^2+t_2^2t_5^2+t_3^2t_4^2)\cdot f_5(q)\\
&+\frac{1}{6}\left[t_1t_6(t_4^3+t_5^3)+t_2t_5(t_4^3+t_6^3)
+t_3t_4(t_5^3+t_6^3)\right]\cdot f_6(q)\\
&+\frac{1}{2}(t_1t_4t_5t_6^2+t_2t_4t_5^2t_6+t_3t_4^2t_5t_6)\cdot f_7(q)\\
&+\frac{1}{4}(t_1t_4^2t_5^2+t_2t_4^2t_6^2+t_3t_5^2t_6^2)\cdot f_8(q)
+\frac{1}{24}(t_1t_6^4+t_2t_5^4+t_3t_4^4)\cdot f_9(q)\\
&+\frac{1}{36}(t_4^3t_5^3+t_4^3t_6^3+t_5^3t_6^3)\cdot f_{10}(q)
+\frac{1}{24}(t_4t_5t_6^4+t_4t_5^4t_6+t_4t_5t_6^4)\cdot f_{11}(q)\\
&+\frac{1}{8}(t_4^2t_5^2t_6^2)\cdot f_{12}(q)
+\frac{1}{720}(t_4^6+t_5^6+t_6^6)\cdot f_{13}(q),
\end{align*}
where $f_i(q), i=0,\dots, 13$ are given by
\begin{align*}
&f_0(q)=\frac{1}{3}\left(\frac{q\frac{d}{dq}a(q)}{1-a(q)^3}\right)^{\frac{1}{2}}
=\frac{\eta(q^9)^3}{\eta(q^3)},
f_1(q)=a(q)f_0(q),\\
&f_2(q)=-\frac{1}{9}\frac{q\frac{d}{dq}f_0}{f_0}+a(q)^2f_0(q)^2,\\
&f_3(q)=f_0(q)^2, 
f_4(q)=a(q)f_0(q)^2,\\
&f_5(q)=-\frac{2}{9}\frac{q\frac{d}{dq}f_0}{f_0}+a(q)^2f_0(q)^2,\\
&f_6(q)=f_0(q)^3,
f_7(q)=a(q)f_0(q)^3, 
f_8(q)=a(q)^2f_0(q)^3,\\
&f_9(q)=a(q)^3f_0(q)^3,
f_{10}(q)=3a(q)f_0(q)^4,
f_{11}(q)=3a(q)^2f_{0}(q)^4,\\
&f_{12}(q)=(2+a(q)^3)f_0(q)^4,
f_{13}(q)=3a(q)(2-a(q)^3)f_{0}(q)^4,
\end{align*}
and 
\[
a(q)=1+\frac{1}{3}\left(\frac{\eta(q)}{\eta(q^9)}\right)^3
=\frac{1}{3}q^{-1}(1+5q^{3}-7q^{6}+3q^{9}+\dots).
\]
%
\end{thm}
\begin{pf}
We can deduce the Theorem from the following uniqueness property of the potential:
\begin{lem}\label{lem:uniquness-E6}
Let $F_0(t_0,\cdots,t_6,t,f_0,\cdots,f_{13})$ 
be a polynomial defined by 
\begin{align*}
&F_0(t_0,\cdots,t_6,t,f_0,\cdots,f_{13})\\
:=&\frac{1}{2}t_0^2 t
+\frac{1}{3}t_0 (t_1 t_6+t_2 t_5+t_3 t_4)
+(t_1 t_2 t_3)\cdot f_0\\
&+\frac{1}{6}(t_1^3+t_2^3+t_3^3)\cdot f_1
+(t_1t_2t_5t_6+t_1t_3t_4t_6+t_2 t_3t_4t_5)\cdot f_2\\
&+\frac{1}{2}(t_1^2t_4t_5+t_2^2t_4t_6+t_3^2t_5t_6)\cdot f_3\\
&+\frac{1}{2}(t_1t_2t_4^2+t_1t_3t_5^2+t_2t_3t_6^2)\cdot f_4
+\frac{1}{4}(t_1^2t_6^2+t_2^2t_5^2+t_3^2t_4^2)\cdot f_5\\
&+\frac{1}{6}\left[t_1t_6(t_4^3+t_5^3)+t_2t_5(t_4^3+t_6^3)
+t_3t_4(t_5^3+t_6^3)\right]\cdot f_6\\
&+\frac{1}{2}(t_1t_4t_5t_6^2+t_2t_4t_5^2t_6+t_3t_4^2t_5t_6)\cdot f_7\\
&+\frac{1}{4}(t_1t_4^2t_5^2+t_2t_4^2t_6^2+t_3t_5^2t_6^2)\cdot f_8
+\frac{1}{24}(t_1t_6^4+t_2t_5^4+t_3t_4^4)\cdot f_9\\
&+\frac{1}{36}(t_4^3t_5^3+t_4^3t_6^3+t_5^3t_6^3)\cdot f_{10}
+\frac{1}{24}(t_4t_5t_6^4+t_4t_5^4t_6+t_4t_5t_6^4)\cdot f_{11}\\
&+\frac{1}{8}(t_4^2t_5^2t_6^2)\cdot f_{12}
+\frac{1}{720}(t_4^6+t_5^6+t_6^6)\cdot f_{13}. 
\end{align*}
\begin{enumerate}
\item For the holomorphic functions 
$f_0(t),\cdots,f_{13}(t)$, 
the holomorphic function \\
$F_0(t_0,\cdots,t_6,t,f_0(t),\cdots,f_{13}(t))$ 
is a potential of an $8$-dimensional Frobenius structure with 
flat coordinates $t_0,t_1,t_2,t_3,t_4,t_5,t_6,t$ 
such that 
the Euler vector field $E$ is given by  
$E=t_0\frac{\p}{\p t_0}+\sum_{k=1}^3\frac{2}{3}t_k\frac{\p}{\p t_k}+\sum_{k=4}^6\frac{1}{3}t_k\frac{\p}{\p t_k}$\\
if and only if there exists $A \in \CC^*$ such that 
\begin{align}
&f_0(t)=A \left(\frac{a(t)'}{1-a(t)^3}\right)^{1/2},
\label{S1202-1}\\ 
&f_1(t)=a(t) f_0(t),\ 
f_2(t)=-\frac{1}{2\cdot 3^2}\left(\frac{a(t)''}{a(t)'}+\frac{a(t)^2 a(t)'}
{1-a(t)^3}\right),\nonumber\\
&f_3(t)=\frac{1}{3^2}\frac{a(t)'}{1-a(t)^3}, \ 
f_4(t)=\frac{1}{3^2}\frac{a(t)a(t)'}{1-a(t)^3},\nonumber\\
&f_5(t)=-\frac{1}{3^2}\left(\frac{a(t)''}{a(t)'}+\frac{2a(t)^2 a(t)'}{1-a(t)^3}\right),\nonumber\\
&f_6(t)=\frac{1}{3^4}A^{-4}f_0(t)^3,
f_7(t)=a(t)f_6(t), 
f_8(t)=a(t)^2f_6(t),\nonumber\\
&f_9(t)=a(t)^3f_6(t),
f_{10}(t)=\frac{1}{3^5}A^{-6}a(t)f_0(t)^4,
f_{11}(t)=a(t) f_{10}(t),\nonumber\\
&f_{12}(t)=\frac{1}{3^6}A^{-6}(2+a(t)^3)f_0(t)^4,
f_{13}(t)=(2-a(t)^3)f_{10}(t),\label{S1202-3}
\end{align}
and
\begin{equation}
\frac{a(t)'''}{a(t)'}-\frac{3}{2}\left(\frac{a(t)''}{a(t)'}\right)^2
=-\frac{1}{2}\frac{8+a(t)^3}{(1-a(t)^3)^2}a(t)\cdot (a(t)')^2, 
\label{S1202-4}
\end{equation}
where $a(t)=f_1(t)/f_0(t)$ and $'=\frac{d}{dt}$. 
\item There exist unique formal power series$:$
\begin{equation}
\widetilde{f}_0(q)=\sum_{n=1}^\infty a_0(n)q^{n},\quad 
\widetilde{f}_i(q)=\sum_{n=0}^\infty a_i(n)q^{n},\quad i=1,\dots, 13, 
\end{equation}
with $a_0(1)=1$ and $a_1(0)=\frac{1}{3}$ such that 
$F_0(t_0,\cdots,t_6,t,\widetilde{f}_0(e^t),\cdots,\widetilde{f}_{13}(e^t))$ 
is the potential of an $8$-dimensional Frobenius structure with 
flat coordinates $t_0$, $t_1,\dots ,t_6,t$, 
and Euler vector field 
$E=t_0\frac{\p}{\p t_0}+\sum_{k=1}^3\frac{2}{3}t_k\frac{\p}{\p t_k}+\sum_{k=4}^6\frac{1}{3}t_k\frac{\p}{\p t_k}$.
\end{enumerate}
\end{lem}
\begin{pf}
The assertion (i) is a direct consequence of WDVV equations 
and discussed already in \cite{E6}. 
For the proof of (ii), we need the following Sub-Lemma.

\begin{sublem}\label{S1129-11}
There exists a unique formal Laurent series
\[
f(q)=\sum_{n=-1}^\infty a_n q^n
\]
satisfying the following conditions$:$
\begin{enumerate}
\item The first coefficient $a_{-1}=\frac{1}{3}$. 
\item $f(q)$ satisfies the following differential equation:
\begin{equation}
\frac{f(q)'''}{f(q)'}-\frac{3}{2}\left(\frac{f(q)''}{f(q)'}\right)^2
=-\frac{1}{2}\frac{8+f(q)^3}{(1-f(q)^3)^2}f(q)\cdot (f(q)')^2, 
\label{S1202-5}
\end{equation}
where $'=q\frac{d}{dq}$. 
\end{enumerate}
\end{sublem}
\begin{pf}
Put \[
S(q):=(1-f(q)^3)^2[f(q)'\cdot f(q)'''-\frac{3}{2}(f(q)'')^2]
+\frac{1}{2}(8+f(q)^3)\cdot f(q) \cdot (f(q)')^4.
\]
Condition (ii) is equivalent to all the coefficients 
of the $q$-expansion of $S(q)$ being zero. 
For the cases of $n \leq 0$, 
the coefficients of $q^{-8+n}$ of $S(q)$ equal to zero. 
For the cases of $n \geq 1$, 
the coefficients of $q^{-8+n}$ of $S(q)$ are of the form
\[
-n^3 a_{-1}^7 a_{n-1}+
\text{ a polynomial in }
a_{-1},\cdots, a_{n-2}.
\] 
Since we have $a_{-1}=1/3$, the coefficients $a_0,a_1,\cdots$ 
are uniquely determined inductively. 
\qed
\end{pf}
We first construct $\widetilde{f}_0(q),\cdots,\widetilde{f}_{13}(q)$. 
Take a formal Laurent series $\widetilde{f}(q)$ as the one which is constructed 
in Sub-Lemma \ref{S1129-11}. 
We take $A \in \CC^*$ such that the formal power series$:$
$A(\frac{q\frac{d}{dq}\widetilde{f}(q)}{1-\widetilde{f}(q)^3})^{1/2}$ 
has an expansion $q+\cdots$. 
Then $A^2$ must be $1/9$. 
We define the following formal power series$:$
\begin{align*}
&\widetilde{f}_0(q):=A(\frac{q\frac{d}{dq}\widetilde{f}(q)}{1-\widetilde{f}(q)^3})^{1/2},\ 
\widetilde{f}_1(q):=\widetilde{f}(q)\widetilde{f}_0(q),\\
&\widetilde{f}_2(q):=-\frac{1}{2\cdot 3^2}
\left(
\frac{(q\frac{d}{dq})^2\widetilde{f}(q)}{q\frac{d}{dq}\widetilde{f}(q)}+
\frac{\widetilde{f}(q)^2q\frac{d}{dq}\widetilde{f}(q)}{1-\widetilde{f}(q)^3}
\right),\cdots
\end{align*}
in a parallel manner as in (\ref{S1202-3}). 
By (i) of this Lemma, we see that 
$\widetilde{f}_i(q)(i=0,\cdots,13)$ satisfy the conditions of (ii). 

We show the uniqueness of $\widetilde{f}_i(q)\ (i=0,\cdots,13)$. 
We assume that $\widehat{f}_i(q)\ (i=0,\cdots,13)$ also satisfy the 
conditions of (ii). 
Put $\widehat{f}(q):=\widehat{f}_1(q)/\widehat{f}_0(q)$. 
By (i) of this Lemma, we see that 
\begin{enumerate}
\item $\widehat{f}(e^t)$ must satisfy the differential 
equation (\ref{S1202-4}). 
\item $\exists \widehat{A} \in \CC^*$ such that 
$$\widehat{f}_0(e^t)=\widehat{A}
\left(
\frac{\frac{d}{dt}\widehat{f}(e^t)}{1-\widehat{f}(e^t)^2}
\right)^{1/2}.$$
\end{enumerate}
From (i), $\widehat{f}(q)$ satisfies (\ref{S1202-5}). 
Since $\widehat{f}(q)$ has the expansion $\frac{1}{3}q^{-1}+\cdots$, 
$\widehat{f}(q)$ must be $\widetilde{f}(q)$ by Sub-Lemma \ref{S1129-11}. 
From (ii) and a comparison of the leading term of $q$-expansions 
of $\widetilde{f}_0(q)$ and $\widehat{f}_0(q)$, we have 
$\widetilde{f}_0(q)=\widehat{f}_0(q)$ and $\widehat{A}^2=A^2$. 
Since $\widehat{f}_i(e^t)\ (i=1,\cdots,13)$ must satisfy 
(\ref{S1202-3}), we have $\widetilde{f}_i(q)=\widehat{f}_i(q)\ (i=1,\cdots,13)$.
Thus we obtain Lemma \ref{lem:uniquness-E6}.
\qed
\end{pf}
Next, we construct the analytic solution to the WDVV equation as follows.
\begin{lem}\label{lem:analytic solution h(q)}
Put 
\begin{equation}
h(q)=1+\frac{1}{3}\left(\frac{\eta(q)}{\eta(q^9)}\right)^3
=\frac{1}{3}q^{-1}+\cdots.
\label{1127-0}
\end{equation}
Then $h(q)$ has the following properties$:$
\begin{enumerate}
\item $h(q)$ satisfies the following differential equation.
\[
\frac{h(q)'''}{h(q)'}-\frac{3}{2}\left(\frac{h(q)''}{h(q)'}\right)^2
=-\frac{1}{2}\frac{8+h(q)^3}{(1-h(q)^3)^2}h(q)\cdot (h(q)')^2, 
\]
where $'=q\frac{d}{dq}$. 
\item $h(q)$ satisfies the following equation$:$	
\begin{equation}
-\frac{1}{64}\frac{h(q)^3(8+h(q)^3)^3}{(1-h(q)^3)^3}=J(q)
\label{1127-2}
\end{equation}
where $J(q)$ is the Laurent series characterized by the conditions that 
\begin{enumerate}
\item $J(q)=\frac{1}{1728}(q^{-3}+744+\cdots)$,
\item $J(\exp(\frac{2\pi\sqrt{-1}\tau}{3}))$ is the elliptic modular 
function on the upper half plane 
$\HH=\{\tau \in \CC\,|\,\mathrm{Im}\tau>0\,\}$. 
\end{enumerate}
\item $h(q)$ has the following expressions$:$
\begin{equation}
h(q)=\omega+\frac{1}{3}\left(\frac{\eta(q\omega^{-2})}{\eta(q^9)}\right)^3
\cdot \exp(\frac{2\pi\sqrt{-1}}{12})
=\omega^2+\frac{1}{3}\left(\frac{\eta(q\omega^{-1})}{\eta(q^9)}\right)^3
\cdot \exp(\frac{2\pi\sqrt{-1}}{24}), 
\label{1127-3}
\end{equation}
where $\omega=\exp(\frac{2\pi\sqrt{-1}}{3})$. 
\end{enumerate}
\end{lem}
\begin{pf}
The uniformization of the Hesse pencil$:$
\[
x_0^3+x_1^3+x_2^3-3ax_0x_1x_2=0
\]
is classically studied and we refer to \cite{S1130-Ohyama2}. 
In \cite{S1130-Ohyama2}, the parameter $a$ is described as a 
holomorphic function $a(\tau)$ on the upper half plane 
$\HH=\{\tau \in \CC\,|\,\mathrm{Im}\tau>0\,\}$ as 
$$
a(\tau)=1+9\left(\frac
{\eta(\exp(2\pi\sqrt{-1}3\tau))}
{\eta(\exp(\frac{2\pi\sqrt{-1}\tau}{3}))}
\right)^3. 
$$
By the modular property of $\eta(\exp(2\pi\sqrt{-1}\tau))$, 
we have 
\begin{equation}\label{eqn:a tau}
a(-\frac{1}{\tau})=h(\exp(\frac{2\pi\sqrt{-1}\tau}{3})).
\end{equation}
Then we can deduce Lemme \ref{lem:analytic solution h(q)} 
from the corresponding results for $a(\tau)$, which 
are classically known and written in \cite{S1130-Ohyama2}. 
\qed
\end{pf}
Finally, we give two important formulas for the function $h(q)$ in Lemma \ref{lem:analytic solution h(q)}:
\begin{lem}\label{lem:3.5}
We have the following equations$:$
\begin{align}
&(1)\ \frac{1}{3^3}\frac{(q\frac{d}{dq}h(q))^6}{(h(q)^3-1)^3}
=\eta(q^3)^{24}.
\label{1127-7}\\
&(2)\ \frac{q\frac{d}{dq}h(q)}{1-h(q)^3}=
3^2\left(\frac{\eta(q^9)^3}{\eta(q^3)}\right)^2.
\label{1127-8}
\end{align}
\end{lem}
\begin{pf}
We have 
\begin{equation}
\frac{1}{2^6\cdot 3^9}\frac{(q\frac{d}{dq}J(q))^6}{J(q)^4(J(q)-1)^3}
=\eta(q^3)^{24},
\label{1127-5}
\end{equation}
because the leading terms of the $q$-expansions coincide and 
if we put $q=\exp(\frac{2\pi\sqrt{-1}\tau}{3})$, 
then both sides are cusp forms of weight 12 
with respect to the $SL(2,\ZZ)$ action and therefore they are uniquely determined 
by the leading terms of the $q$-expansions. 

By (\ref{1127-2}) and (\ref{1127-5}), we have (\ref{1127-7}). 

We could easily check that 
\begin{equation}
\exp(\frac{2\pi\sqrt{-1}}{24})\eta(q)\eta(q\omega^{-1})
\eta(q\omega^{-2})\eta(q^9)
=(\eta(q^3))^4.
\label{1127-4}
\end{equation}
By (\ref{1127-0}), (\ref{1127-3}), (\ref{1127-4}), we have 
\begin{equation}
h(q)^3-1=\frac{1}{3^3}\left(\frac{\eta(q^3)}{\eta(q^9)}\right)^{12}.
\label{1127-6}
\end{equation}

By (\ref{1127-7}), (\ref{1127-6}) and the comparison of 
the leading terms of $q$-expansions, we have (\ref{1127-8}). 
\qed
\end{pf}

It is easy to show that, by our choice of basis $\gamma_0,\dots, \gamma_7$ of the orbifold cohomology group 
$H^*_{orb}(\PP^1_{3,3,3},\QQ)$ and their dual coordinates $t_0,\dots, t_7$ in the beginning of this section,
the Gromov--Witten potential is of the form in Lemma \ref{lem:uniquness-E6} 
except for the condition $a_0(1)=1$. 
Indeed, we can choose elements $\gamma_1,\gamma_6\in H^*_{orb}(\PP^1_{3,3,3},\CC)$ contained in 
the basis which will correspond to coordinates $t_1,t_6$ such that 
$\gamma_1\circ \gamma_1=\gamma_6$ and 
$\int_{\PP^1_{3,3,3}}\gamma_1\cup \gamma_6=\frac{1}{3}$ where $\circ$ denotes 
the orbifold cohomology ring structure on $H^*_{orb}(\PP^1_{3,3,3},\CC)$. 
This gives us $a_1(0)=\frac{1}{3}$.
The condition $a_0(1)=1$ follows from the fact that 
the Gromov--Witten invariant $a_0(1)$ counts the number of  
morphisms from $\PP^1_{3,3,3}$ to $\PP^1_{3,3,3}$ of degree one,
which is exactly the identity map. Hence, we have $a_0(1)=1$.
Now, the statement in Theorem \ref{thm:genus zero E6} follows from the uniquness of the potential.
\qed
\end{pf}
Now, we consider the Frobenius structure 
on the base space of the universal unfolding of simple elliptic singularity of 
type $\widetilde{E}_6: W_{\widetilde{E}_6}(x_1,x_2,x_3):=x_1^3+x_2^3+x_3^3-3ax_1x_2x_3$. 
It is easily obtained once we fix a primitive form (see \cite{st:1} for example).
It is proven by K.~Saito in \cite{sa:1} that there exists a primitive form 
for $W_{\widetilde{E}_6}(x_1,x_2,x_3)=x_1^3+x_2^3+x_3^3-3ax_1x_2x_3$ 
and it is given by choosing a cycle in the corresponding elliptic curve 
$\{W_{\widetilde{E}_6}(x_1,x_2,x_3)=0\}\subset \PP^2$.
Denote by $M_{\widetilde{E}_6,\infty}$ the Frobenius manifold 
with the choice of the primitive form associated to the cycle in the elliptic curve 
which vanishes when the parameter $a$ goes to infinity. 
In view of (i) of Lemma \ref{lem:uniquness-E6},
we only have to calculate the holomorphic function $a(t)$ in order to 
describe the potential for $M_{\widetilde{E}_6,\infty}$.
However, it is also easy to see from the result in \cite{sa:1} that 
we can choose the uniformization parameter $\tau/3$ 
as the flat coordinate $t$ for our choice of primitive form
and hence we have $a(\tau)=h(\exp(\frac{2\pi\sqrt{-1}\tau}{3}))$ 
as in the equation \eqref{eqn:a tau}.
By rescaling other flat coordinates suitably, 
it is possible to set $A=1/3$ (in the notation of (i) of Lemma \ref{lem:uniquness-E6}). 
Therefore, we can apply the uniqueness of the potential, (ii) of Lemma \ref{lem:uniquness-E6}, 
and hence we obtain an isomorphism $M_{\PP^1_{3,3,3}}\simeq M_{\widetilde{E}_6,\infty}$ 
as Frobenius manifolds. 
On the other hand, for the elliptic root system of type $E_6^{(1,1)}$
(\cite{extendedII}), 
the domain $\mathbb E_{E_6^{(1,1)}}$ and 
the elliptic Weyl group $W_{E_6^{(1,1)}}$ are defined and 
the quotient space $M_{E_6^{(1,1)}}:={\mathbb E}_{E_6^{(1,1)}}//W_{E_6^{(1,1)}}
\simeq 
\{z \in \CC\,|\,\text{Re}z<0\,\}\times \CC^7$
has a Frobenius manifold structure isomorphic 
to $M_{\widetilde{E_6},\infty}$ (\cite{S1202-Saito}, 
\cite{extendedII}, \cite{SatakeFrobenius}).
To summarize, we obtain the following
\begin{thm}
We have isomorphisms of Frobenius manifolds
\[
M_{\PP^1_{3,3,3}}\simeq M_{\widetilde{E}_6,\infty}\simeq M_{E_6^{(1,1)}}.
\]
\qed
\end{thm}
\subsection{Genus one potential}
We shall also give the genus one Gromov--Witten potential.
\begin{thm}
The genus one Gromov--Witten potential $F_1^{\PP^1_{3,3,3}}$ of $\PP^1_{3,3,3}$ 
is given as
\begin{equation}
F_1^{\PP^1_{3,3,3}}=-\frac{1}{3}\log(\eta(q^3)).
\end{equation}
\end{thm}
\begin{pf}
The proof is similar to the one for $F_1^{\PP^1_{2,2,2,2}}$.
It is easy to see that the genus one Gromov--Witten potential 
$F_1^{\PP^1_{3,3,3}}$ is an element of $\QQ[[q]]$ 
since the Euler vector field is given by 
$E=t_0\frac{\p}{\p t_0}+\sum_{k=1}^3\frac{2}{3}t_k\frac{\p}{\p t_k}+\sum_{k=4}^6\frac{1}{3}t_k\frac{\p}{\p t_k}$.
Therefore, we only have to consider the (orbifold) stable maps 
with one marked point from 
smooth elliptic curves to $\PP^1_{3,3,3}=[{\bf E}/(\ZZ/3\ZZ)]$,
which factor through ${\bf E}$ by definition.
Hence, we have that 
\[
q\frac{d}{d q}F_1^{\PP^1_{2,2,2,2}}=
-\frac{1}{24}+\sum_{n=1}^{\infty}\sigma(n)q^{3n}.
\]
One may also obtain the statement by Dubrovin--Zhang's Virasoro constraint \cite{dz:1}.
Indeed, Proposition 4 in \cite{dz:1} gives us the equation
\[
q\frac{d}{d q}F_1^{\PP^1_{3,3,3}}=\frac{3}{4}f_2(q)+\frac{3}{8}f_5(q)=
-\frac{1}{12}q\frac{d}{dq}\log\left(q\frac{d}{dq}h(q)\right)
-\frac{1}{8}\frac{q\frac{d}{dq}h(q)\cdot h(q)^2}{1-h(q)^3}.
\]
By the equation (\ref{1127-7}) in Lemma \ref{lem:3.5}, 
we have $-\frac{1}{12}q\frac{d}{dq}\log(q\frac{d}{dq}h(q))
-\frac{1}{8}\frac{q\frac{d}{dq}h(q)\cdot h(q)^2}{1-h(q)^3}=-\frac{1}{3}
q\frac{d}{dq}\log(\eta(q^3))$.
\qed
\end{pf}
Strachan \cite{strachan:1} calculates the G-function 
for the Frobenius structure on the universal unfolding of simple elliptic singularities of 
type $\widetilde{E}_6, \widetilde{E}_7, \widetilde{E}_8$
with the choice of the primitive form ``at $a=0$". 
If we use the primitive form ``at $a= \infty$" instead, 
then G-functions for $\widetilde{E}_6$, $\widetilde{E}_7$ and $\widetilde{E}_8$
can be obtained as 
$-\frac{1}{3}\log(\eta(q^3))$, $-\frac{1}{4}\log(\eta(q^4))$ and $-\frac{1}{4}\log(\eta(q^4))$
respectively.
This is consistent with our calculation of Gromov--Witten invariants and mirror symmetry.


\begin{thebibliography}{99}
{\small 
\bibitem{agv:1}
D.~Abramovich, T.~Graber and A.~Vistoli,
	{\it Gromov--Witten theory of Deligne--Muford stacks},
	Amer. J. Math. {\bf 130} (2008), no. 5, 1337--1398.
\bibitem{cr:1}
W.~Chen and Y.~Ruan,
	{\it Orbifold Gromov--Witten Theory},
	Orbifolds in mathematics and physics (Madison, WI, 2001), 25--85, Contemp. Math., 
	{\bf 310}, Amer. Math. Soc., Providence, RI, 2002.
\bibitem{d:1}
B.~Dubrovin,
{\it Geometry of 2d topological field theories},
Integrable systems and quantum groups
(Montecatini Terme, 1993), Lecture Notes in Math., vol. 1620, Springer, Berlin, 1996, pp. 120--348.
\bibitem{dz:1}
B.~Dubrovin and Y.~Zhang,
{\it Frobenius manifolds and Virasoro constraints},
Selecta Math. (N.S.) {\bf 5} (1999), no. 4, 423--466. 
\bibitem{dz:2}
B.~Dubrovin and Y.~Zhang,
{\it Bihamiltonian hierarchies in 2D topological field theory at
one-loop approximation}. 
Commun. Math. Phys. {\bf 198} (1998), 311 -- 361.
\bibitem{et:1}
W.~Ebeling and A.~Takahashi
{\it Strange duality of weighted homogeneous polynomials},
Compositio Math., {\bf 147}, 1413 -- 1433, (2011).
\bibitem{fjr:1}
H.~Fan, T.~Jarvis and Y.~Ruan,
{\it The Witten equation, mirror symmetry and quantum singularity theory}, 
arXiv:0712.4021.
\bibitem{ks:1}
M.~Krawitz, Y.~Shen,
{\it Landau-Ginzburg/Calabi-Yau Correspondence of all Genera
for Elliptic Orbifold $\PP^1$}. 
arXiv:1106.6270.
%
%
%
%
\bibitem{mr:1}
T.~Milanov, Y.~Ruan, 
{\it Gromov-Witten theory of elliptic orbifold $\PP^1$ and quasimodular forms}.
arXiv:1106.2321.
\bibitem{Mumford}
D.~Mumford, 
\textit{Tata Lectures on Theta I}, 
Progress in Math. {\bf 28}, 
Birkh\"auser (1983). 
\bibitem{ny:1}
M.~Noumi and Y.~Yamada,
{\it Notes on the flat structures associated with simple and simply elliptic
singularities}, 
in ``Integrable Systems and Algebraic Geometry" (eds. M.-H. Saito, Y. Shimuzu, K. Ueno),
Proceedings of the Taniguchi Symposium 1997, 373--383, World-Scientific, 1998.

\bibitem{S1130-Ohyama1}
Y.~Ohyama, 
\textit{Differential relations of theta functions}. 
Osaka J. Math. 32 (1995), no. 2, 431--450. 
%
%
\bibitem{S1130-Ohyama2}
Y.~Ohyama, 
\textit{Differential equations for modular forms of level three}. 
Funkcial. Ekvac. 44 (2001), no. 3, 377--389.
%
%
%
%
\bibitem{r:1}
P.~Rossi,
{\it Gromov-Witten theory of orbicurves, the space of tri-polynomials and Symplectic Field Theory of Seifert fibrations},
arXiv:0808.2626. 

\bibitem{sa:1}
K.~Saito, 
{\it Primitive forms for a universal unfolding of a function with an isolated critical point}. 
J. Fac. Sci. Univ. Tokyo Sect. IA Math. {\bf 28} (1982), no. 3, 775--792. 
%
%
\bibitem{S1202-Saito}
K.~Saito, 
\textit{Period mapping associated to a primitive form}, 
Publ. RIMS, Kyoto Univ. {\bf 19} (1983) 1231--1264.
%
%
\bibitem{extendedII}
K.~Saito, 
\textit{Extended Affine Root System II}, 
Publ. RIMS. Kyoto Univ. {\bf 26} (1990) 15--78.
\bibitem{sa:eta}
K.~Saito, {\it Duality for regular systems of weights}, 
Asian J. Math. {\bf 2} (1998) 983--1047.
\bibitem{st:1}
K.~Saito and A.~Takahashi,
{\it From Primitive Forms to Frobenius manifolds}, 
Proceedings of Symnposia in Pure Mathematics, {\bf 78} (2008) 31--48.
\bibitem{D4}
I.~Satake, 
\textit{Flat Structure and the Prepotential 
for the Elliptic Root System of Type $D^{(1,1)}_4$}, 
Topological field theory, primitive forms 
and related topics (ed. by M. Kashiwara, et al.), 
Progress in Math. {\bf 160}, Birkh\"auser (1998), 427--452. 
%
%
\bibitem{SatakeFrobenius}
I.~Satake, 
\textit{Frobenius manifolds for elliptic root systems}, 
Osaka J. Math. {\bf 47} (2010) 301--330.
%
%
\bibitem{strachan:1}
I.~Strachan,
{\it Simple elliptic singularities: A note on their G-function},
arXiv:1004.2140.
\bibitem{t:1}
A.~Takahashi,
{\it Weighted projective lines associated to regular systems of weights of dual type},
Adv. Stud. Pure Math. {\bf 59} (2010), 371--388.
\bibitem{u:1}
K.~Ueda,
{\it Homological mirror symmetry and simple elliptic singularities},
arXiv:math/0604361.
\bibitem{E6}
E.~Verlinde, N.P.~Warner, 
\textit{Topological Landau-Ginzburg matter at $c=3$}, 
Physics Letters B {\bf 269} (1991) 96--102. 
}
\end{thebibliography}
\end{document}